\documentclass{article}
\usepackage{amsmath,amssymb}
\newcommand{\Z}{\ensuremath{\mathbf Z}}
\newcommand{\N}{\ensuremath{ \mathbf N }}
\newcommand{\R}{\ensuremath{ \mathbf R }}
\newtheorem{theorem}{Theorem}
\newcommand{\bt}{\begin{theorem}}
\newcommand{\et}{\end{theorem}}
\newtheorem{lemma}{Lemma}
\newcommand{\bl}{\begin{lemma}}
\newcommand{\el}{\end{lemma}}
\newcommand{\pf}{{\bf Proof}.\ }
\newcommand{\be}{\begin{eqnarray}}
\newcommand{\ee}{\end{eqnarray}}
\newcommand{\beq}{\begin{equation}}
\newcommand{\eeq}{\end{equation}}
\newcommand{\benum}{\begin{enumerate}}
\newcommand{\eenum}{\end{enumerate}}
\newcommand{\bal}{\begin{align*}}
\newcommand{\eal}{\end{align*}}
\newcommand{\ba}{\begin{array}}
\newcommand{\ea}{\end{array}}
\newcommand{\eop}{$\square$\vspace{.3cm}}
\DeclareMathOperator{\card}{card}

\date{}

\begin{document}
\title{Every function is the representation function\\
of an additive basis for the integers\footnote{2000 Mathematics
Subject Classification:  11B13, 11B34, 11B05.
Key words and phrases.  Additive bases, sumsets, representation functions, 
density, Erd\H os-Tur\' an conjecture, Sidon set.}}
\author{Melvyn B. Nathanson\thanks{This work was supported
in part by grants from the NSA Mathematical Sciences Program
and the PSC-CUNY Research Award Program.}\\
Department of Mathematics\\
Lehman College (CUNY)\\
Bronx, New York 10468\\
Email: nathansn@alpha.lehman.cuny.edu}
\maketitle

\begin{abstract}
Let $A$ be a set of integers.
For every integer $n$, let $r_{A,h}(n)$ denote
the number of representations of $n$ in the form
$ n =  a_1 + a_2 + \cdots +  a_h,$ where $a_1, a_2, \ldots,a_h \in A$ 
and $a_1 \leq a_2 \leq \cdots \leq a_h.$
The function 
\[
r_{A,h}: \Z \rightarrow \N_0\cup\{\infty\}
\]
is the {\em representation function of order $h$ for $A$.}
The set $A$ is called an {\em asymptotic basis of order $h$} if 
$r_{A,h}^{-1}(0)$ is finite, that is, if every integer 
with at most a finite number of exceptions can be represented 
as the sum of exactly $h$ not necessarily distinct elements of $A$.
It is proved that every function is a representation function, that is,
if $f: \Z \rightarrow \N_0\cup\{\infty\}$ is any function such that 
$f^{-1}(0)$ is finite, 
then there exists a set $A$ of integers such that $f(n) = r_{A,h}(n)$
for all $n \in \Z$.
Moreover, the set $A$ can be arbitrarily sparse in the sense that,
if $\varphi(x) \geq 0$ for $x \geq 0$ and $\varphi(x) \rightarrow \infty$, 
then there exists a set $A$ with $f(n) = r_{A,h}(n)$ 
and $\card\left(\{a\in A : |a| \leq x\}\right) < \varphi(x)$ for all $x$. 

It is an open problem to construct dense sets of integers 
with a prescribed representation function.  
Other open problems are also discussed.
\end{abstract}

\section{Additive bases and the Erd\H os-Tur\' an conjecture}
Let $\N, \N_0$, and \Z\ denote the positive integers, 
nonnegative integers, and integers, respectively.
Let $A$ be a set of integers.
For every positive integer $h$, we define the {\em sumset}
\[
hA = \{a_1+\cdots + a_h : a_i \in A \text{ for all } i = 1,\ldots, h\}.
\]
We denote by $r_{A,h}(n)$ the number
of representations of $n$ in the form $n = a_1 + a_2 + \cdots + a_h,$ 
where $a_1, a_2, \ldots,a_h \in A$ and $a_1 \leq a_2 \leq \cdots \leq a_h.$
The function $r_{A,h}$ is called the {\em representation function of order $h$} 
of the set $A$.

We denote the cardinality of the set $A$ by $\card(A)$.
If $A$ is a finite set of integers, we denote the maximum element 
of $A$ by $\max(A).$
For any integer $t$ and set $A \subseteq \Z$, we define the {\em translate}
\[
t+A = \{t+a : a \in A\}.
\]

In this paper we consider additive bases for the set of all integers.
The set $A$ of integers is called a {\em basis of order $h$ for \Z}
if every integer can be represented as the sum of $h$
not necessarily distinct elements of $A$.
The set $A$ of integers is called an 
{\em asymptotic basis of order $h$ for \Z}
if every integer with at most a finite number of exceptions 
can be represented as the sum of $h$
not necessarily distinct elements of $A$.
Equivalently, the set $A$ is an asymptotic basis of order $h$ 
if the representation function 
$r_{A,h}:\Z \rightarrow \N_0 \cup \{\infty\}$ satisfies the condition 
\[
\card\left(r_{A,h}^{-1}(0)\right) < \infty.
\]

For any set $X$, let $\mathcal{F}_0(X)$ denote the set of all functions 
\[
f:X \rightarrow \N_0 \cup\{\infty\}
\]
such that
\[
\card\left(f^{-1}(0)\right) < \infty.
\]
We ask:  Which functions in $\mathcal{F}_0(\Z)$
are representation functions of asymptotic bases for the integers?  
This question has a remarkably simple and surprising answer.
In the case $h=1$ we observe that $f \in \mathcal{F}_0(\Z)$ 
is a representation function if and only if $f(n) = 1$  
for all integers $n \not\in f^{-1}(0).$
For $h \geq 2$ we shall prove that {\em every}\, function in $\mathcal{F}_0(\Z)$ 
is a representation function.  Indeed, if $f \in \mathcal{F}_0(\Z)$ and $h \geq 2,$ 
then there exist infinitely many sets $A$ such that
$f(n) = r_{A,h}(n)$ for every $n \in \Z$.
Moreover, we shall prove that we can construct arbitrarily sparse 
asymptotic bases $A$ with this property.  
Nathanson~\cite{nath03a} previously proved this theorem for
$h=2$ and the function $f(n) = 1$ for all $n \in \Z.$

This result about asymptotic bases for the integers contrasts sharply with the
case of the nonnegative integers. 
The set $A$ of nonnegative integers is called 
an {\em asymptotic basis of order $h$ for $\N_0$}
if every sufficiently large integer can be represented as the sum of $h$
not necessarily distinct elements of $A$.
Very little is known about the class of representation functions
of asymptotic bases for $\N_0$.
However, if $f\in \mathcal{F}_0(\N_0),$ then Nathanson~\cite{nath78c}
proved that there exists at most one set $A$ such that $r_{A,h}(n) = f(n).$

Many of the results that have been proved about asymptotic bases for $\N_0$ are negative.
For example, Dirac~\cite{dira51} showed that the representation function 
of an asymptotic basis of order 2 cannot be eventually constant.
Erd\H os and Fuchs~\cite{erdo-fuch56} proved that the average value of a representation
function of order 2 cannot even be approximately constant, in the sense that,
for every infinite set $A$ of nonnegative integers and every real number $c > 0$,
\[
\sum_{n\leq N} r_{A,2}(n) \neq cN + o\left( N^{1/4}\log^{-1/2}N\right).
\]
Erd\H os and Tur\' an~\cite{erdo-tura41} conjectured that 
if $A$ is an asymptotic basis of order $h$ for the nonnegative integers,
then the representation function $r_{A,h}(n)$ must be unbounded,
that is, 
\[
\limsup_{n\rightarrow\infty} r_{A,h}(n) = \infty.
\]
This famous unsolved problem in additive number theory is only 
a special case of the general problem of classifying the representation functions 
of asymptotic bases of finite order for the nonnegative integers.

\section{Two lemmas}

We use the following notation.
For sets $A$ and $B$ of integers and for any integer $t$,
we define the {\em sumset}
\[
A+B = \{a+b: a \in A, b \in B\},
\]
the {\em translation}
\[
A + t = \{a+t: a \in A\},
\]
and the {\em difference set}
\[
A - B = \{a - b: a \in A, b \in B\}.
\]
For every nonnegative integer $h$ we define the {\em $h$-fold sumset $hA$}
by induction:
\begin{align*}
0A & = \{0\},  \\
hA & = A + (h-1)A = \{a_1 + a_2 + \cdots + a_h : a_1, a_2,\ldots, a_h \in A\}.
\end{align*}
The {\em counting function} for the set $A$ is
\[
A(y,x) = \card\left( \{a\in A: y \leq a \leq x\}\right).
\]
In particular, $A(-x,x)$ counts the number of integers $a\in A$
with $|a| \leq x.$

Let $[x]$ denote the integer part of the real number $x$.

\bl           \label{Zreps:lemma:uk}
Let $f:\Z \rightarrow \N_0\cup\{\infty\}$ be a function such that $f^{-1}(0)$ is finite.  
Let $\Delta$ denote the cardinality of the set $f^{-1}(0).$
Then there exists a sequence $U = \{u_k\}_{k=1}^{\infty}$
of integers such that, for every $n \in \Z$ and $k \in \N$,
\[
f(n) = \card\left( \{ k\geq 1 : u_k = n \}\right)
\]
and
\[
|u_k| \leq \left[\frac{k+\Delta}{2}\right].
\]
\el

\pf
Every positive integer $m$ can be written uniquely in the form 
\[
m = s^2+s+1+r, 
\]
where $s$ is a nonnegative integer and $|r|\leq s.$
We construct the sequence 
\begin{align*}
V & = \{0,-1,0,1,-2,-1,0,1,2,-3,-2,-1,0,1,2,3,\ldots\} \\
& = \{v_m\}_{m=1}^{\infty},
\end{align*}
where 
\[
v_{s^2+s+1+r} = r \qquad\mbox{for $|r|\leq s.$}
\]
For every nonnegative integer $k$, the first occurrence of $-k$ 
in this sequence is $v_{k^2+1} =  -k,$
and the first occurrence of $k$ in this sequence is $v_{(k+1)^2} =  k.$

The sequence $U$ will be the unique subsequence of $V$ constructed as follows.
Let $n \in \Z.$  If $f(n) = \infty,$ then $U$ will contain 
the terms $v_{s^2+s+1+n}$ for every $s \geq |n|$.
If $f(n) = {\ell} <\infty,$ then $U$ will contain 
the ${\ell}$ terms $v_{s^2+s+1+n}$ for $s = |n|, |n|+1,\ldots,|n|+{\ell}-1$ in the subsequence $U$,
but not the terms $v_{s^2+s+1+n}$ for $s \geq |n| + {\ell}.$
Let $m_1 < m_2 < m_3 < \cdots$ be the strictly increasing sequence of positive integers
such that $\{v_{m_k}\}_{k=1}^{\infty}$ is the resulting subsequence of $V$.
Let $U = \{u_k\}_{k=1}^{\infty}$, where $u_k = v_{m_k}.$
Then 
\[
f(n) = \card\left( \{ k \geq 1 : u_k = n \}\right).
\]

Let $\card\left(f^{-1}(0)\right) = \Delta.$
The sequence $U$ also has the following property:
If $|u_k| = n,$ then for every integer $m \not\in f^{-1}(0)$ 
with $|m| < n$ there is a positive integer $j < k$ with $u_j = m$.  
It follows that
\[
\{0,1,-1,2,-2,\ldots, n-1, -(n-1)\}\setminus f^{-1}(0) \subseteq \{u_1,u_2,\ldots,u_{k-1}\},
\]
and so
\[
k-1 \geq 2(n-1)+1 - \Delta.
\]
This implies that
\[
|u_k| = n \leq \frac{k+\Delta}{2}.
\]
Since $u_k$ is an integer, we have
\[
|u_k| \leq \left[\frac{k+\Delta}{2}\right].
\]
This completes the proof.
\eop

Lemma~\ref{Zreps:lemma:uk} is best possible in the sense that
for every nonnegative integer $\Delta$ there is a function
$f:\Z \rightarrow \N_0\cup\{\infty\}$ with $\card\left(f^{-1}(0)\right) = \Delta$ 
and a sequence $U = \{u_k\}_{k=1}^{\infty}$
of integers such that
\beq               \label{Zreps:U}
|u_k| = \left[\frac{k+\Delta}{2}\right]  \qquad \text{for all $k \geq 1.$}
\eeq
For example, if $\Delta = 2\delta + 1$ is odd, define the function $f$ by
\[
f(n) = \left\{\ba{ll}
0 & \text{if $|n| \leq \delta$}  \\
1 & \text{if $|n| \geq \delta + 1$} 
\ea
\right.
\]
and the sequence $U$ by
\begin{align*}
u_{2i-1} & = \delta + i, \\
u_{2i} & = -(\delta + i)
\end{align*}
for all $i \geq 1.$

If $\Delta = 2\delta$ is even, define $f$ by
\[
f(n) = \left\{\ba{ll}
0 & \text{if $-\delta \leq n \leq \delta-1$}  \\
1 & \text{if $n \geq \delta$ or $n \leq -\delta-1$} 
\ea
\right.
\]
and the sequence $U$ by $u_1 = \delta$ and
\begin{align*}
u_{2i} & = \delta + i, \\
u_{2i+1} & = -(\delta + i)
\end{align*}
for all $i \geq 1.$
In both cases the sequence $U$ satisfies~(\ref{Zreps:U}).

The set $A$ is called a {\em Sidon set of order $h$} 
if $r_{A,h}(n) = 0 \text{ or } 1$ for every integer $n$.  
If $A$ is a Sidon set of order $h$, then $A$ is a Sidon set
of order $j$ for all $j = 1,2,\ldots, h.$

\bl           \label{Zreps:lemma:sidon}
Let $A$ be a finite Sidon set of order $h$ and $d = \max(\{|a| : a \in A\}).$
If $|c| > (2h-1)d,$ then $A \cup \{ c\}$ is also a Sidon set of order $h$.
\el

\pf
Let $n \in h\left( A \cup \{c\}\right).$  Suppose that
\[
n = a_1+\cdots + a_j + (h-j)c = a'_1 + \cdots + a'_{{\ell}} + (h-{\ell})c,
\]
where 
\[
0 \leq j \leq {\ell} \leq h,
\]
\[
a_1,\ldots,a_j,a'_1,\ldots,a'_{\ell} \in A,
\]
and
\[
a_1 \leq \cdots \leq a_j \qquad\text{and}\qquad a'_1\leq \cdots\leq a'_{\ell}.
\]
If $j < {\ell},$ then
\begin{align*}
|c| & \leq |({\ell}-j)c|  \\
& = \left| a'_1 + \cdots + a'_{\ell} - (a_1+\cdots + a_j) \right|  \\ 
& \leq ({\ell}+j)d  \\
& \leq (2h-1)d  \\
& < |c|,
\end{align*}
which is absurd.
Therefore, $j={\ell}$ and $a_1+\cdots + a_j = a'_1 + \cdots + a'_j$.
Since $A$ is a Sidon set of order $j$, 
it follows that $a_i = a'_i$ for all $i = 1,\ldots,j.$ 
Consequently, $A\cup \{c\}$ is a Sidon set of order $h$.
\eop

\section{Construction of asymptotic bases}
We can now construct asymptotic bases of order $h$ for the integers 
with arbitrary representation functions.

\bt\label{Zreps:theorem:main}
Let $f:\Z \rightarrow \N_0\cup\{\infty\}$ be a function
such that the set $f^{-1}(0)$ is finite.
Let $\varphi:\N_0 \rightarrow\R$ be a nonnegative function such that 
$\lim_{x\rightarrow\infty} \varphi(x) = \infty.$
For every $h \geq 2$ there exist infinitely many asymptotic bases $A$ 
of order $h$ for the integers such that
\[
r_{A,h}(n) = f(n) \qquad\text{for all $n \in \Z$,}
\]
and
\[
A(-x,x) \leq \varphi(x)
\]
for all $x \geq 0$.
\et

\pf
By Lemma~\ref{Zreps:lemma:uk}, there is a sequence 
$U = \{u_k\}_{k=1}^{\infty}$ of integers such that
\[
f(n) = \card\left(\{ k \geq 1 : u_k = n\}\right)
\]
for every integer $n$.

Let $h \geq 2.$  We shall construct an ascending sequence of finite sets
$A_1 \subseteq A_2 \subseteq A_3 \subseteq \cdots$
such that, for all positive integers $k$ and for all integers $n,$
\benum
\item[(i)]
\[
r_{A_k,h}(n) \leq f(n) ,
\]
\item[(ii)]
\[
r_{A_{k},h}(n) \geq \card\left(\{i : 1 \leq i \leq k \mbox{ and } u_i = n\}\right),
\]
\item[(iii)]
\[
\card(A_k) \leq 2k,
\]
\item[(iv)]
\[
\text{$A_k$ is a Sidon set of order $h-1$.}
\]
\eenum
Conditions~(i) and~(ii) imply that the infinite set
\[
A = \bigcup_{k=1}^{\infty} A_k
\]
is an asymptotic basis of order $h$ for the integers 
such that $r_{A,h}(n) = f(n)$ for all $n \in \Z$.

We construct the sets $A_k$ by induction.
Since the set $f^{-1}(0)$ is finite, there exists an integer $d_0$ such that 
$f(n) \geq 1$ for all integers $n$ with $|n| \geq d_0.$
If $u_1 \geq 0,$ choose a positive integer $c_1 > 2hd_0$.
If $u_1 < 0,$ choose a negative integer $c_1 < -2hd_0$.
Then
\[
|c_1| > 2hd_0.
\]
Let 
\[
A_1 = \{-c_1, (h-1)c_1+u_1  \}.
\]  
The sumset $hA_1$ is the finite arithmetic progression
\begin{align*}
hA_1 & = \{ -hc_1 + (hc_1+u_1)i : i = 0,1,\ldots,h\} \\
& = \{ -hc_1, u_1, hc_1+2u_1,2hc_1+3u_1,\ldots,  (h-1)hc_1+hu_1\}.
\end{align*}
Then $|n| \geq hc_1 > d_0$ for all $n \in hA_1\setminus\{u_1\},$ and so
\[
r_{A_1,h}(n) = 1 \leq f(n)
\]
for all $n \in hA_1$.
It follows that $r_{A_1,h}(n) \leq f(n)$ for all $n \in \Z$.
The set $A_1$ is a Sidon set of order $h$, hence also a Sidon set of order $h-1.$
Thus, the set $A_1$ satisfies conditions~(i)--(iv).

We assume that for some integer $k \geq 2$ we have
constructed a set $A_{k-1}$ satisfying conditions~(i) -- (iv).
If
\[
r_{A_{k-1}}(n) \geq \card\left( \{i : 1 \leq i \leq k \mbox{ and } u_i = n\}\right)
\]
for all $n\in \Z,$ then the set $A_{k} = A_{k-1}$
satisfies conditions~(i) -- (iv).
Otherwise, 
\[
r_{A_{k-1}}(u_{k}) 
= \card\left(\{i : 1 \leq i \leq k \mbox{ and } u_i = u_{k}\}\right) - 1 < f(u_{k}).
\]
We shall construct a Sidon set $A_{k}$ of order $h-1$ 
such that 
\[
\card(A_{k}) = \card(A_{k-1})+2
\]
and 
\beq          \label{Zreps:rep}
r_{A_{k},h}(n) = \left\{
\begin{array}{ll}
r_{A_{k-1},h}(n)+1 & \text{if $n = u_{k}$}\\
r_{A_{k-1},h}(n) & \text{if $n \in hA_{k-1} \setminus \{u_{k}\}$}\\
1 & \text{if $n \in hA_{k} \setminus \left( hA_{k-1} \cup \{u_{k}\}\right).$}
\end{array}
\right.
\eeq

Define the integer
\beq           \label{Zreps:dk}
d_{k-1} = \max\left(\left\{ |a| : a\in A_{k-1} \cup \{u_{k}\}\right\}\right).
\eeq
Then
\[
A_{k-1} \subseteq [-d_{k-1},d_{k-1}].
\]
If $u_{k} \geq 0,$ choose a positive integer $c_{k}$ such that
$c_{k} > 2hd_{k-1}.$
If $u_{k} < 0,$ choose a negative integer $c_{k}$ 
such that $c_{k} < -2hd_{k-1}.$
Then 
\beq            \label{Zreps:ck}
|c_k| > 2hd_{k-1}.
\eeq

Let 
\[
A_{k} = A_{k-1} \cup \{-c_k,(h-1)c_k+u_{k}\}.
\]
Then
\[
\card(A_{k}) = \card(A_{k-1})+2 \leq 2k.
\]

We shall assume that $u_k \geq 0$, hence $c_k > 0.$
(The argument in the case $u_k < 0$ is similar.)
We decompose the sumset $hA_k$ as follows:
\[
hA_{k} 
= \bigcup_{\substack{r+i+j=h\\ r,i,j\geq 0}} 
\left( r(h-1)c_k + ru_{k} -ic_k + jA_{k-1}  \right)
= \bigcup_{r=0}^h B_r,
\]
where
\[
B_r = r(h-1)c_k + ru_{k} + \bigcup_{i=0}^{h-r}\left( -ic_k + (h-r-i)A_{k-1}  \right).
\]
If $n \in B_r$, then there exist integers $i \in \{ 0,1,\ldots,h-r\}$ 
and $y \in (h-r-i)A_{k-1}$ such that
\[
n = r(h-1)c_k + ru_k - ic_k + y.
\]
Since
\[
|y| \leq (h-r-i)d_{k-1},
\]
it follows that
\beq            \label{Zreps:nlow}
n \geq r(h-1)c_k + ru_k - ic_k - (h-r-i)d_{k-1}
\eeq
and
\[
n \leq r(h-1)c_k + ru_k -ic_k + (h-r-i)d_{k-1}.
\]

Let $m \in B_{r-1}$ and $n \in B_r$ for some $r \in \{1,\ldots,h\}.$ 
There exist nonnegative integers $i \leq h-r$ and
$j \leq h-r+1$ such that 
\begin{align*}
n-m 
& \geq \left(r(h-1)c_k+ru_k - ic_k - (h-r-i)d_{k-1} \right) \\
& \qquad - \left((r-1)(h-1)c_k+(r-1)u_k -jc_k + (h-r+1-j)d_{k-1} \right)\\
& = (h-1+j-i)c_k+u_k  - (2h-2r-i-j+1)d_{k-1} \\
& \geq (h-1-i)c_k - 2hd_{k-1}.
\end{align*}
If $r \geq 2,$ then $i \leq h-2$ and inequality~(\ref{Zreps:ck}) implies that
\[
n-m \geq c_k - 2hd_{k-1} > 0.
\]
Therefore, if $m \in B_{r-1}$ and $n \in B_r$ for some $r \in \{2,\ldots,h\},$ 
then $m < n$. 

In the case $r=1$ we have $m \in B_0$ and $n \in B_1$.  
If $i \leq h-2$, then~(\ref{Zreps:ck}) implies that
\[
n-m  \geq (h-1-i)c_k - 2hd_{k-1} \geq c_k - 2hd_{k-1} > 0
\]
and~(\ref{Zreps:nlow}) implies that
\[
n \geq (h-1-i)c_k + u_k - (h-1-i)d_{k-1} > c_k - hd_{k-1} > d_0.
\]
If $r=1$ and $i = h-1$, then $n = u_k$.  
Therefore, if $m \in B_0$ and $n\in B_1$, then $m < n$ unless $m = n = u_k.$
It follows that the sets $B_0 ,B_1\setminus \{u_k\} ,B_2,\ldots, B_h$ are pairwise disjoint.

Let $n \in B_r$ for some $r \geq 1$.  Suppose that $0 \leq i \leq j \leq h-r,$ 
and that 
\[
n = r(h-1)c_k+ru_k -ic_k + y  \qquad \text{for some $y \in (h-r-i)A_{k-1}$} 
\]
and
\[
n = r(h-1)c_k+ru_k - jc_k + z       \qquad \text{for some $z \in (h-r-j)A_{k-1}$.}
\]
Subtracting these equations, we obtain
\[
z-y = (j-i)c_k.
\]
Recall that $|a| \leq d_{k-1}$ for all $a \in A_{k-1}.$
If $i < j,$ then
\begin{align*}
c_k & \leq (j-i)c_k = z-y \\
& \leq |y|+|z| \leq (2h-2r-i-j)d_{k-1} \\
& < 2hd_{k-1}  < c_k,
\end{align*}
which is impossible.  Therefore, $i = j$ and $y = z$.  
Since $0 \leq h-r-i \leq h-1$ and $A_{k-1}$ is a Sidon set of order $h-1$, it follows that 
\[
r_{A_{k-1},h-r-i}(y) = 1 
\]
and so
\[
r_{A_k,h}(n) = 1 \leq f(n) \qquad
\text{for all }  n \in \left( B_1\setminus\{u_k\}\right) \cup \bigcup_{r=2}^{h}B_r.
\]

Next we consider the set
\[
B_0 = hA_{k-1} \cup \bigcup_{i=1}^h \left( -ic_k + (h-i)A_{k-1}\right).
\]
For $i = 1,\ldots,h$, we have 
\[
c_k > 2hd_{k-1} \geq (2h-2i+1)d_{k-1}
\] 
and so
\begin{align*}
\max\left( -ic_k + (h-i)A_{k-1}\right) & \leq -ic_k+(h-i)d_{k-1} \\
& < -(i-1)c_k-(h-i+1)d_{k-1} \\
& \leq \min\left( -(i-1)c_k + (h-i+1)A_{k-1}  \right).
\end{align*}
Therefore, the sets $-ic_k + (h-i)A_{k-1} $ are pairwise disjoint
for $i = 0,1,\ldots,h.$  
In particular, if $n \in B_0\setminus hA_{k-1},$ then
\[
n \leq \max\left( -c_k + (h-1)A_{k-1}\right) \leq -c_k+(h-1)d_{k-1} < -d_{k-1} \leq -d_0
\]
and $f(n) \geq 1.$
Since $A_{k-1}$ is a Sidon set of order $h-1$,
it follows that 
\[
r_{A_k,h}(n) = 1 \leq f(n) \qquad\text{for all }
n \in \bigcup_{i=1}^h \left(-ic_k + (h-i)A_k\right)= B_0\setminus hA_{k-1}.
\]
It follows from~(\ref{Zreps:dk}) that for any $n \in B_0\setminus hA_{k-1}$ we have 
\[
n < -d_{k-1} \leq u_k,
\]
and so $u_k \not\in B_0\setminus hA_{k-1}.$ 
Therefore,
\[
r_{A_k,h}(u_k) = r_{A_{k-1},h}(u_k)+1, 
\]
and the representation function $r_{A_k,h}$ satisfies
the three requirements of~(\ref{Zreps:rep}).

We shall prove that 
\[
A_{k} = A_{k-1} \cup \{-c_k,(h-1)c_k+u_{k}\}.
\]
is a Sidon set of order $h-1$.
Since $A_{k-1}$ is a Sidon set of order $h-1$ 
with $d_{k-1} \geq \max\{ |a|: a \in A_{k-1}\}$,
and since 
\[
c_k > 2hd_{k-1} > (2(h-1)-1)d_{k-1},
\]
Lemma~\ref{Zreps:lemma:sidon} implies
that $A_{k-1} \cup \{-c_k\}$ is a Sidon set of order $h-1$.

Let $n \in (h-1)A_k.$
Suppose that
\begin{align*}
n & = r(h-1)c_k + ru_k - ic_k + x \\
& = s(h-1)c_k + su_k - jc_k + y, 
\end{align*}
where 
\[
0 \leq r \leq s \leq h-1,
\] 
\[
0 \leq i \leq h-1-r,
\]
\[
0 \leq j \leq h-1-s,
\]
\[
x \in (h-1-r-i)A_{k-1},
\]
and
\[
y \in (h-1-s-j)A_{k-1}.
\]
Then
\[
|x| \leq (h-1-r-i)d_{k-1}
\]
and
\[
|y| \leq (h-1-s-j)d_{k-1}.
\]
If $r < s,$ then $j \leq h-2$ and 
\begin{align*}
(h-1)c_k & \leq (s-r)(h-1)c_k + (s-r)u_k \\
 & = (j-i)c_k + x-y \\
 & \leq (j-i)c_k + (2h-2-r-s-i-j)d_{k-1} \\
& \leq (h-2)c_k + 2h d_{k-1} \\
& < (h-1)c_k,
\end{align*}
which is absurd.  Therefore, $r=s$ and
\[
- ic_k + x = - jc_k + y \in (h-1-r)\left( A_k \cup \{-c_k\} \right).
\]
Since $A_k \cup \{-c_k\}$ is a Sidon set of order $h-1$, it follows that $i= j$
and that $x$ has a unique representation as the sum of $h-1-r-i$ elements of $A_k.$
Thus, $A_k$ is a Sidon set of order $h-1$.

The set $A_k$ satisfies conditions~(i) -- (iv).  It follows by induction that there 
exists an infinite increasing sequence $A_1 \subseteq A_2 \subseteq \cdots $
of finite sets with these properties, and that $A = \cup_{k=1}^{\infty} A_k$
is an asymptotic basis of order $h$ with representation function $r_{A,h}(n) = f(n)$
for all $n \in \Z.$

Let $A_0 = \emptyset,$ and let $K$ be the set of all positive integers $k$
such that $A_k \neq A_{k-1}.$  Then
\[
A = \cup_{k \in K} A_k = \cup_{k \in K} \{ -c_k, (h-1)c_k \}.
\]
For $k \in K,$ 
the only constraints on the choice of the number $c_k$ 
in the construction of the set $A_{k}$ 
were the sign of $c_k$ and the growth condition
\[
|c_k| > 2hd_{k-1}  \qquad\text{for all integers $k \in K$.}
\]
We shall prove that we can construct the asymptotic basis $A$ 
with counting function $A(-x,x) \leq \varphi(x)$ for all $x \geq 0.$
Since $\varphi(x) \rightarrow \infty$ as $x \rightarrow \infty$,
for every integer $k \geq 0$ there exists an integer $w_k$
such that 
\[
\varphi(x) \geq 2k \qquad\text{for all $x \geq w_k.$}
\]
We now impose the following additional constraint:
Choose $c_k$ such that 
\[
|c_k| \geq w_k \qquad\text{for all integers $k \in K$.}
\]
Then
\[
A_1(-x,x) = 0 \leq \varphi(x) \qquad\text{for $0 \leq x < |c_1|$}
\]
and
\[
A_1(-x,x) \leq 2 \leq \varphi(x) \qquad\text{for $ x \geq |c_1| \geq w_1.$}
\]
Suppose that $k \geq 2$ and the set $A_{k-1}$ satisfies $A_{k-1}(-x,x) \leq \varphi(x)$
for all $x \geq 0.$ 
Since
\[
A_k \cap (-|c_k|,|c_k|) = A_{k-1} \cap (-|c_k|,|c_k|),
\]
it follows that
\[
A_k(-x,x) =  A_{k-1}(-x,x) \leq \varphi(x) \qquad\text{for $0 \leq x < |c_k|$}
\]
and
\[
A_k(-x,x) \leq 2k \leq \varphi(x) 
\qquad\text{for $ x \geq |c_k| \geq w_k.$}
\]
This proves by induction that $A_k(-x,x) \leq \varphi(x)$
for all $k$ and $x$.  Since $\lim_{k\rightarrow\infty} |c_k| = \infty,$ 
it follows that for any nonnegative integer $x$ we can choose $c_k$ so that
$|c_k| > x$ and 
\[
A(-x,x) = A_k(-x,x) \leq \varphi(x).
\]
For every integer $k \in K$ we had infinitely many choices for the integer $c_k$
to use in the construction of the set $A_k,$ and so there are 
infinitely many asymptotic bases $A$ with the property that
$r_A(n) = f(n)$ for all $n \in \Z.$
This completes the proof.
\eop

\section{Sums of pairwise distinct integers}

Let $A$ be a set of integers and $h$ a positive integer.
We define the sumset $h\wedge A$ as the set consisting of all sums 
of $h$ pairwise distinct elements of $A$, 
and the {\em restricted representation function}
\[
\hat{r}_{A,h}: \Z \rightarrow \N_0 \cup \{\infty\}
\]
by
\[
\hat{r}_{A,h}(n) 
= \card\left( \left\{ \{a_1,\ldots,a_h\} \subseteq A: a_1+\cdots + a_h = n 
\text{ and } a_1 < \cdots < a_h \right\} \right).
\]
The set $A$ of integers is called a {\em restricted asymptotic basis of order $h$}
if $h\wedge A$ contains all but finitely many integers, or, equivalently, 
if $\hat{r}_{A,h}^{-1}(0)$ is a finite subset of \Z.

The following theorem can also be proved by the method used to prove 
Theorem~\ref{Zreps:theorem:main}.

\bt\label{Zreps:theorem:distinct}
Let $f:\Z \rightarrow \N_0\cup\{\infty\}$ be a function
such that $f^{-1}(0)$ is a finite set of integers.
Let $\varphi:\N_0 \rightarrow \R$ be a nonnegative function such that 
$\lim_{x\rightarrow\infty} \varphi(x) = \infty.$
For every $h \geq 2$ there exist infinitely many sets $A$ of integers such that
\[
\hat{r}_{A,h}(n) = f(n) \qquad\text{for all $n \in \Z$}
\]
and
\[ 
A(-x,x) \leq \varphi(x)
\]
for all $x \geq 0$.
\et

\section{Open problems}
Let $X$ be an abelian semigroup, written additively, 
and let $A$ be a subset of $X$.
We define the $h$-fold sumset $hA$ as the set consisting 
of all sums of $h$ not necessarily distinct elements of $A$.
The set $A$ is called an {\em asymptotic basis of order $h$ for $X$ }
if the sumset $hA$ consists of all but at most finitely many elements of $X$.
We also define the $h$-fold {\em restricted sumset} $h\wedge A$ as the set consisting 
of all sums of $h$ pairwise distinct elements of $A$.
The set $A$ is called a {\em restricted asymptotic basis of order $h$ for $X$ }
if the restricted sumset $h\wedge A$ consists of all 
but at most finitely many elements of $X$.
The classical problems of additive number theory concern the semigroups $\N_0$ and $\Z$.

There are four different representation functions that we can
associate to every subset $A$ of $X$ and every positive integer $h$.
Let $(a_1,\ldots,a_h)$ and $(a'_1,\ldots,a'_h)$ be $h$-tuples of elements of $X$.
We call these $h$-tuples {\em equivalent} if there is a permutation $\sigma$ 
of the set $\{1,\ldots,h\}$ such that $a'_{\sigma(i)} = a_i$ for all $i = 1,\ldots, h.$
For every $x \in X,$ let $r_{A,h}(x)$ denote the number of equivalence classes 
of $h$-tuples $(a_1,\ldots,a_h)$ of elements of $A$ such that
$a_1+\cdots + a_h = x.$  The function $r_{A,h}$ 
is called the {\em unordered representation function} of $A$. 
This is the function that we studied in this paper.
The set $A$ is an asymptotic basis of order $h$ if
$r_{A,h}^{-1}(0)$ is a finite subset of $X$.

Let $R_{A,h}(x)$ denote the number of $h$-tuples $(a_1,\ldots,a_h)$ 
of elements of $A$ such that
$a_1+\cdots + a_h = x.$  The function $R_{A,h}$ 
is called the {\em ordered representation function} of $A$. 

Let $\hat{r}_{A,h}(x)$ denote the number of equivalence classes 
of $h$-tuples $(a_1,\ldots,a_h)$ of pairwise distinct elements of $A$ 
such that $a_1+\cdots + a_h = x$, and let $\hat{R}_{A,h}(x)$ 
denote the number of $h$-tuples $(a_1,\ldots,a_h)$ 
of pairwise distinct elements of $A$ such that $a_1+\cdots + a_h = x.$  
These functions are called the 
{\em unordered restricted representation function} of $A$
and the {\em ordered restricted representation function} of $A$, respectively.
The two restricted representation functions are essentially identical, since 
$\hat{R}_{A,h}(x) = h!\hat{r}_{A,h}(x)$ for all $x \in X.$

In the discussion below, we use only the unordered representation function 
$r_{A,h},$ but each of the problems can be reformulated in terms of the
other representation functions.

For every countable abelian semigroup $X$, 
let $\mathcal{F}(X)$ denote the set of all functions
$f: X \rightarrow \N_0 \cup \{\infty\}$, and let 
$\mathcal{F}_0(X)$ denote the set of all functions
$f: X \rightarrow \N_0 \cup \{\infty\}$ such that $f^{-1}(0)$ 
is a finite subset of $X$.
Let $\mathcal{F}_c(X)$ denote the set of all functions
$f: X \rightarrow \N_0 \cup \{\infty\}$ such that $f^{-1}(0)$ 
is a cofinite subset of $X$,
that is, $f(x) \neq 0$ for only finitely many $x \in X,$ or, equivalently,
\[
\card\left( f^{-1}(\N \cup \{\infty\}) \right) < \infty.
\]

Let $\mathcal{R}(X,h)$ denote the set of all $h$-fold representation functions 
of subsets $A$ of $X$.
If $r_{A,h}$ is the representation function of an asymptotic basis $A$ 
of order $h$ for $X$, then $r_{A,h}^{-1}(0)$ is a finite subset of $X$, and so 
$r_{A,h} \in \mathcal{F}_0(X)$.
Let $\mathcal{R}_0(X,h)$ denote the set of all $h$-fold representation functions 
of asymptotic bases $A$ of order $h$ for $X$.
Let $\mathcal{R}_c(X,h)$ denote the set of all $h$-fold representation functions 
of finite subsets of $X$.
We have
\[
\mathcal{R}(X,h) \subseteq \mathcal{F}(X),
\]
\[
\mathcal{R}_0(X,h) \subseteq \mathcal{F}_0(X),
\]
and
\[
\mathcal{R}_c(X,h) \subseteq \mathcal{F}_c(X),
\]

In the case $h=1,$ we have, for every set $A \subseteq X$,
\[
r_{A,1}(x) = \left\{
\begin{array}{ll}
1 & \text{if $x \in A$,}\\
0 & \text{if $x \not\in A$,}
\end{array}
\right.
\]
and so 
\[
\mathcal{R}(X,1)= \{f:X\rightarrow \{0,1\}\},
\]
\[
\mathcal{R}_0(X,1)= \{f:X\rightarrow \{0,1\}: \card\left(f^{-1}(0)\right) < \infty \},
\]
and
\[
\mathcal{R}_c(X,1)= \{f:X\rightarrow \{0,1\}: \card\left(
f^{-1}(\N \cup \{\infty\}) \right) < \infty \}.
\]
In this paper we proved that
\[
\mathcal{R}_0(\Z,h) = \mathcal{F}_0(\Z)
\qquad\text{for all $h \geq 2.$}
\]
Nathanson~\cite{nath04a} has extended this result to certain 
countably infinite groups and semigroups.
Let $G$ be any countably infinite abelian group
such that $\{2g:g\in G\}$ is infinite. 
For the unordered restricted representation function $\hat{r}_{A,2},$
we have
\[
\mathcal{R}_0(G,2) = \mathcal{F}_0(G).
\]
More generally, let $S$ is any countable abelian semigroup 
such that for every $s \in S$ there exist $s', s'' \in S$ with $s = s' + s''.$
In the abelian semigroup $X = S\oplus G,$ we have
\[
\mathcal{R}_0(X,2) = \mathcal{F}_0(X).
\]
If $\{12g:g\in G\}$ is infinite, then $\mathcal{R}_0(X,2) = \mathcal{F}_0(X)$
for the unordered representation function $r_{A,2}.$

The following problems are open for all $h \geq 2$:
\benum
\item
Determine $\mathcal{R}_0(\N_0,h).$  
Equivalently, describe the representation functions of additive bases
for the nonnegative integers.
This is a major unsolved problem in additive number theory, 
of which the Erd\H os-Tur\' an conjecture is only a special case.

\item
Determine $\mathcal{R}(\Z,h).$  
In this paper we computed $\mathcal{R}_0(\Z,h),$ 
the set of representation functions of additive bases
for the integers, 
but it is not known under what conditions a function 
$f:\Z\rightarrow \N_0\cup \{\infty\}$ with $f^{-1}(0)$ infinite 
is the representation function of a subset $A$ of $X$.
It can be proved that if $f^{-1}(0)$ is infinite but sufficiently sparse,
then $f \in \mathcal{R}(\Z,h).$

\item
Determine $\mathcal{R}(\N_0,h).$  Is there a simple list of necessary and sufficient 
conditions for a function $f:\N_0 \rightarrow \N_0$ to be the 
representation function of some set of nonnegative integers?

\item
Determine $\mathcal{R}_c(\Z,h).$  
Describe the representation functions of finite sets of integers.
If $A$ is a finite set of integers and $t$ is an integer,
then for the translated set $t+A$ we have
\[
r_{t+A,h}(n) = r_{A,h}(n-ht)
\]
for all integers $n$.  
This implies that if $f(n) \in \mathcal{R}_c(\Z,h),$ 
then $f(n - ht) \in \mathcal{R}_c(\Z,h)$ 
for every integer $b$, so it suffices to consider only finite sets $A$
of nonnegative integers with $0 \in A,$ 
and functions $f \in \mathcal{F}_c(\N_0,h)$ with $f(0) = 1.$ 

\item
Determine $\mathcal{R}_0(G,2)$, $\mathcal{R}(G,2)$,
and $\mathcal{R}_c(G,2)$ for the infinite abelian group
$G = \oplus_{i=1}^{\infty}\Z/2\Z.$
Note that $\{2g : g \in G\} = \{ 0 \}$ for this group.

\item
Determine $\mathcal{R}_0(G,h)$ and $\mathcal{R}(G,h)$, where $G$ 
is an arbitrary countably infinite abelian group and $h \geq 2.$

\item
There is a class of problems of the following type.
Do there exist integers $h$ and $k$ with $2 \leq h < k$
such that
\[
\mathcal{R}(\Z,h) \neq \mathcal{R}(\Z,k)?
\]
We can easily find sets of integers to show that that 
$\mathcal{R}_0(\N_0,h) \neq \mathcal{R}_0(\N_0,k).$
For example, let $A = \N$ be the set of all positive integers, and let $h \geq 1$.
Then $r_{\N,h}(0) = 0$ and $r_{\N,h}(h) = 1.$
If $B$ is any set of nonnegative integers and $k > h$, 
then either $r_{B,k}(0) = 1$ or $r_{B,k}(h) = 0$,
and so $r_{\N,h} \not\in \mathcal{R}_0(\N_0,k).$
Is it true that
\[
\mathcal{R}_0(\N_0,h) \cap \mathcal{R}_0(\N_0,k) = \emptyset
\]
for all $h \neq k$?

\item
By Theorem~\ref{Zreps:theorem:main}, for every $h \geq 2$ 
and every function $f \in \mathcal{F}_0(\Z)$, there exist arbitrarily
sparse sets $A$ of integers such that $r_{A,h}(n) = f(n)$ for all $n.$
It is an open problem to determine how dense the sets $A$ can be.
For example, in the special case $h = 2$ and $f(n) = 1,$
Nathanson~\cite{nath03a} proved that there exists a set $A$ such that 
$r_{A,2}(n) = 1$ for all $n,$ and $\log x \ll A(-x,x) \ll \log x$.  
For an arbitrary representation function $f \in \mathcal{F}_0(\Z)$, 
Nathanson~\cite{nath03f} constructed an asymptotic basis of order $h$ with 
$A(-x,x) \gg x^{1/(2h-1)}.$ 
In the case $h = 2,$ Cilleruelo and Nathanson~\cite{cill-nath04}
improved this to $A(-x,x) \gg x^{\sqrt{2}-1}$.

\eenum

\providecommand{\bysame}{\leavevmode\hbox to3em{\hrulefill}\thinspace}
\providecommand{\MR}{\relax\ifhmode\unskip\space\fi MR }
\providecommand{\MRhref}[2]{%
  \href{http://www.ams.org/mathscinet-getitem?mr=#1}{#2}
}
\providecommand{\href}[2]{#2}

\end{document}